\documentclass[12pt]{amsart} 

\usepackage{amsmath, amsthm, amscd, amsfonts, amssymb, enumerate,verbatim, newlfont, calc, graphicx, color}
\usepackage[bookmarksnumbered, colorlinks, plainpages]{hyperref}
\usepackage{tikz} 
\usepackage{doi}
\usepackage{cancel}
\newtheorem{theorem}{Theorem}[section] 
 
\newtheorem{fact}[theorem]{Fact} 
 
\newtheorem{lemma}[theorem]{Lemma}
\newtheorem{proposition}[theorem]{Proposition}

\newtheorem{corollary}[theorem]{Corollary}

\newtheorem{example}[theorem]{Example}

\usepackage{mathtools}

\usepackage{pstricks}
\usepackage{epsfig}
\usepackage{pst-grad}
\usepackage{pst-plot}
\usepackage{subfig}

\begin{document}

\author[M. Nasernejad,  V. Crispin   Qui$\mathrm{\tilde{n}}$onez, and   W. Hochst\"attler]{Mehrdad ~Nasernejad$^{1,*}$,  Veronica Crispin Qui$\mathrm{\tilde{n}}$onez$^{2}$,  Winfried  Hochst\"attler$^{3}$}

\title[Normally torsion-free square-free monomial ideals] {On the normally torsion-freeness of square-free monomial ideals}



\thanks{$^*$Corresponding author}

\thanks{E-mail addresses:   m$\_$nasernejad@yahoo.com,  veronica.crispin@math.uu.se, and winfried.hochstaettler@fernuni-hagen.de}  
\maketitle

\begin{center}
{\it
$^{1}$Univ. Artois, UR 2462, Laboratoire de Math\'{e}matique de \\ Lens (LML),  F-62300 Lens, France \\

 $^{2}$Department of Mathematics, Uppsala University,  S-751 06, \\ Uppsala, Sweden\\

$^{3}$FernUniversit\"{a}t in Hagen, Fakult\"{a}t f\"{u}r Mathematik und   \\  Informatik,  58084 Hagen,  Germany

 }
\end{center}

\vspace{0.4cm}


\begin{abstract}
Let $I\subset R=K[x_1, \ldots, x_n]$  be a   square-free monomial ideal,  $\mathfrak{q}$ be  a prime monomial ideal in $R$,    $h$ be a square-free monomial in $R$ with  $\mathrm{supp}(h) \cap  (\mathrm{supp}(\mathfrak{q}) \cup \mathrm{supp}(I))=\emptyset$, and $L:=I\cap (\mathfrak{q}, h)$. In this paper, 
we first focus on the associated primes of powers of $L$ and explore the normally torsion-freeness of $L$. We also  give an   application on a combinatorial result. Next, we  study when a square-free monomial ideal is minimally not normally torsion-free. Particularly, we introduce  a class of square-free monomial ideals, which are minimally not normally torsion-free.
\end{abstract}


\vspace{0.4cm}


\section{Preliminaries}

Let  $R$ be  a commutative Noetherian ring and $I$ be   an ideal of $R$. A prime ideal $\mathfrak{p}\subset  R$ is an {\it associated prime} 
    of $I$ if there exists an element $v$ in $R$ such that $\mathfrak{p}=(I:_R v)$. The  {\it set of associated primes} of $I$, denoted by  $\mathrm{Ass}_R(R/I)$, is the set of all prime ideals associated to  $I$.  A well-known result of  Brodmann \cite{BR} showed that the sequence $\{\mathrm{Ass}_R(R/I^k)\}_{k \geq 1}$ of associated prime ideals is stationary  for large $k$. That is, there exists a positive integer $k_0$ such that $\mathrm{Ass}_R(R/I^k)=\mathrm{Ass}_R(R/I^{k_0})$ for all $k\geq k_0$.  The  minimal such $k_0$ is called the
{\it index of stability}   of  $I$ and $\mathrm{Ass}_R(R/I^{k_0})$ is called the {\it stable set } 
 of associated prime ideals of  $I$, which is denoted by $\mathrm{Ass}^{\infty }(I).$

 An ideal $I\subset R$ is said to be  {\it normally torsion-free} if, for all $k\geq 1$, we have  $\mathrm{Ass}_R(R/I^k)\subseteq \mathrm{Ass}_R(R/I)$, 
  see \cite[Definition 4.3.28]{V1}.  It is well-known that a finite simple undirected graph is bipartite if and only if its edge ideal is normally torsion-free, if and only if its cover ideal is normally torsion-free, consult \cite{GRV, SVV}  for more information. In addition, in \cite[Corollary 4.6]{HRV},  it has been shown  that every  transversal polymatroidal ideal  is normally torsion-free.

  It has been proved  in \cite[Theorem 3.2]{KHN1} that   the Alexander dual of path ideals generated by all paths of length 2 in rooted trees are normally
torsion-free. Also, in \cite[Theorem  3.3]{N3}, the authors showed  that  the Alexander dual of the monomial
ideal generated by the paths of maximal lengths in  a rooted starlike tree is normally torsion-free. 
Furthermore, in  \cite[Lemma 5.15]{NQKR}, it has been presented a  characterization of  all normally torsion-free $t$-spread principal Borel ideals. 
More recently, in  \cite[Theorem 4.6]{MNQ}, the authors proved that the edge ideal of any  $\mathbf{t}$-spread $d$-partite hypergraph   
 $K^{\textbf{t}}_{V}$ is normally torsion-free. 
  
   However, normally torsion-freeness  has been explored in the literature, but  little is known for the   normally torsion-free monomial ideals which are not square-free, for instance refer to  \cite{OL}.

In this work, we  continue  studying   the concept of normally torsion-freeness of  square-free monomial ideals. 
For this purpose,  we start to  prove the first main result of this paper in Lemma  \ref{NTF1}.  
As an application of Lemma \ref{NTF1}, we re-prove  that   dominating  ideals of trees are normally torsion-free, see Theorem \ref{DI-TREES},  which   has been shown in \cite[Theorem 2.3]{NQ} by combinatorial tools.   Next, we  establish the second main result of this text in  Lemma  \ref{Ass-maximal1}. 
 
In broadly speaking,   normally torsion-free square-free monomial ideals are closely related to Mengerian hypergraphs.  A hypergraph  is called 
\textit{Mengerian} if it satisfies a certain min-max equation, which is known as the Mengerian property in hypergraph theory or has the max-flow min-cut property in integer programming. In fact,  it is well-known that if $\mathcal{C}$ is  a clutter and $I=I(\mathcal{C})$  its edge ideal, then $I$ is  normally torsion-free if and only if  $I^{k}=I^{(k)}$ for all $k\geq 1$, if and only if $\mathcal{C}$ is Mengerian, if and only if $\mathcal{C}$  has the max-flow min-cut  (MFMC for short)  property,  see  \cite[Theorem 14.3.6]{V1}. On the other hand, one of the famous and long-standing questions in this subject  is 
\textit{Conforti-Cornu\'e{j}ols conjecture}, which says that a hypergraph has the packing property if and only if it is Mengerian, see 
\cite[Conjecture 14.3.19]{V1}. Along this conjecture, the following definition can be arisen:

  A clutter $\mathcal{C}$ is {\it minimally non-MFMC} if it does not have the MFMC property
but all its proper minors do. Equivalently,  $I(\mathcal{C})$  is  {\it  minimally not normally torsion-free} if  $I(\mathcal{C})$  is not normally torsion-free but all its proper minors are  normally torsion-free.  

In the  third main result of this paper,  we  concentrate  on the above  definition,  see Proposition  \ref{MNNT}. In particular, we introduce  a class of square-free monomial ideals originating from graph theory, which are minimally not normally torsion-free,     refer to  Example \ref{Exam.MNNT}.  
 
 \bigskip
Throughout this paper, $R=K[x_1, \ldots, x_n]$ stands for a polynomial  ring over a field $K$ in the $n$  indeterminates $x_1, \ldots, x_n$. 
 Moreover, $\mathcal{G}(I)$ denotes the unique minimal set of monomial generators of a monomial ideal $I\subset R$.  
The {\em support} of a monomial $u\in R$, denoted by $\mathrm{supp}(u)$, is the set of variables that divide $u$. Moreover, for a monomial ideal $I$, we set $\mathrm{supp}(I)=\bigcup_{u \in \mathcal{G}(I)}\mathrm{supp}(u)$. Furthermore,  when there is no confusion about the underlying ring, we will denote the set of associated primes of $I$ simply by $\mathrm{Ass}(R/I)$ or $\mathrm{Ass}(I)$. Also, any  necessary definitions related to graph theory  can be found in  \cite{West}.

The following facts will be used frequently in the rest of this paper.  In order to  ease of reference and more readability, we cite them here.

\begin{fact}\label{proposition 3.3}  (\cite[Proposition 3.3]{NQ})
Suppose that  $I$ and  $J$  are   two  normally torsion-free square-free monomial ideals in a polynomial ring  $R=K[x_1, \ldots, x_n]$ over a field $K$ such that  $\mathrm{supp}(I) \cap \mathrm{supp}(J)=\emptyset$.  Then $I\cap J=IJ$ is normally torsion-free.
\end{fact}

\begin{fact} \label{lemma 2.17}(\cite[Lemma 2.17]{NQBM})
Let $I$ be a normally torsion-free square-free  monomial ideal in a polynomial ring $R=K[x_{1},\ldots ,x_{n}]$ with $\mathcal{G}(I) \subset R$. Then the ideal $L:=IS\cap (x_{n}, x_{n+1}, x_{n+2}, \ldots, x_m)\subset  S=R[x_{n+1}, x_{n+2}, \ldots, x_m]$  is normally torsion-free.
\end{fact}

Note that  we  denote by $V^*(I)$ the set of monomial prime ideals containing $I$.

 \begin{fact}\label{lemma 4.6(viii)}(\cite[Lemma 4.6(viii)]{RNA})
Let $I$ be a monomial ideal   in $R=K[x_1, \ldots, x_n]$, and $\mathfrak{p}$ be a monomial prime ideal of $R$.  If $\mathfrak{p} \in V^*(I)$, then  $\mathrm{Ass}_{R(\mathfrak{p})}(R(\mathfrak{p})/I(\mathfrak{p}))=\{\mathfrak{q}~:~ \mathfrak{q}\in \mathrm{Ass}_{R}(R/I)~\mathrm{and}~ \mathfrak{q} \subseteq \mathfrak{p}\}.$
\end{fact}

 \begin{fact}\label{lemma 4.10(viii)}(\cite[Lemma 4.10(viii)]{RNA})
Let $I$ be a monomial ideal   in   $R=K[x_1,\ldots, x_n]$ over a field $K$, $\mathfrak{m}=(x_1, \ldots, x_n)$, and 
$1\leq j \leq n$.  If $I\subseteq \mathfrak{m}\setminus\{x_j\}$, then $$\mathrm{Ass}_{R/x_j}((R/x_j)/(I/x_j))=\{\mathfrak{q}~:~ \mathfrak{q}\in \mathrm{Ass}_{R}(R/I)~\mathrm{and}~ x_j\notin\mathfrak{q}\}.$$
\end{fact}

\begin{fact} \label{lemma 3.12}(\cite[Lemma 3.12]{SN})
Let  $I$ be a monomial ideal in a polynomial ring $R=K[x_1, \ldots, x_n]$ with 
 $G(I)=\{u_1, \ldots, u_m\}$, and $h=x_{j_1}^{b_1}\cdots x_{j_s}^{b_s}$ with $j_1, \ldots, j_s \in \{1, \ldots, n\}$ be a monomial in $R$. Then  $I$ is normally torsion-free  if and only if $hI$ is normally torsion-free. 
\end{fact}

\begin{fact} \label{theorem 3.15}(\cite[Theorem 3.15]{SN})
Let $I$ be a  monomial ideal in a polynomial ring  $R=K[x_1,\ldots, x_n]$, and $\mathfrak{p}\in V^*(I)$. If $I$ is normally torsion-free, 
then   $I(\mathfrak{p})$  is so.
\end{fact}

\begin{fact} \label{theorem 3.19}(\cite[Theorem 3.19]{SN})
Let $I$ be  a    monomial ideal in  $R=K[x_1, \ldots, x_n]$, and     $1\leq j \leq n$.  If $I$  is normally torsion-free, then    $I/x_j$ is so.
\end{fact}

\begin{fact}\label{theorem 3.3} (\cite[Theorem  3.3]{SNQ})
 Let $I_1\subset R_1=K[x_1, \ldots, x_n]$ and  
$I_2\subset R_2=K[y_1, \ldots, y_m]$ be two monomial ideals in disjoint sets of variables. Let 
$I=I_1R+I_2R\subset R=K[x_1, \ldots, x_n, y_1, \ldots, y_m].$
Then $\mathfrak{p}\in \mathrm{Ass}(R/I)$ if and only if 
$\mathfrak{p}=\mathfrak{p}_1R + \mathfrak{p}_2R$, where 
$ \mathfrak{p}_1\in \mathrm{Ass}(R_1/I_1)$ 
and $ \mathfrak{p}_2\in \mathrm{Ass}(R_2/I_2)$.
\end{fact}

\begin{fact} \label{corollary 3.5} (\cite[Corollary 3.5]{SNQ})
 Let $I\subset  R=K[x_1, \ldots, x_n]$ be a square-free  monomial ideal,  $\mathfrak{m}=(x_1, \ldots, x_n)$,   
 and  $\{u_1, \ldots, u_{\beta_1(I)}\}$  be a maximal independent set of minimal generators of $I$ such that  $\mathfrak{m}\setminus x_i \notin \mathrm{Ass}(R/(I\setminus x_i)^t)$ for all   $x_i\mid \prod_{i=1}^{\beta_1(I)}u_i$ and some positive integer $t$, where $I\setminus x_i$ denotes the  deletion of $I$ at $x_i$. 
 If $\mathfrak{m}\in \mathrm{Ass}(R/I^t)$, then $t\geq \beta_1(I)+1$. 
 \end{fact}

\begin{fact} \label{proposition 4.3.29} (\cite[Proposition 4.3.29]{V1})
 Let $I$ be an ideal of a ring $R$. If $I$ has no embedded
primes, then $I$ is normally torsion-free if and only if $I^n = I^{(n)}$ for all $n \geq 1$. 
\end{fact}

\begin{fact} \label{exercise 6.1.25} (\cite[Exercise 6.1.25]{V1})
 If $\mathfrak{q}_1, \ldots, \mathfrak{q}_r$  are primary monomial ideals of $R$ with non-comparable
radicals and $I$ is an ideal such that $I=\mathfrak{q}_1 \cap \cdots \cap  \mathfrak{q}_r$, then
$I^{(n)}=\mathfrak{q}^n_1 \cap \cdots \cap  \mathfrak{q}^n_r$. 
\end{fact}


\section{Some results on  normally torsion-free square-free monomial ideals }

In this section, we present some results on  normally torsion-free square-free monomial ideals. 
To show  Proposition \ref{PRODUCT}, we employ   Lemma \ref{PRODUCT1}. We start with  proving  Lemma \ref{PRODUCT1} in the following result.

\begin{lemma}\label{PRODUCT1}
Let $I\subset R=K[x_1, \ldots, x_r]$  be a  square-free monomial ideal  and  $\mathfrak{q}$ be  a prime monomial ideal in $R$ such that $\bigcap_{\mathfrak{p}\in \mathrm{Ass}(I)}\mathfrak{p} \cap \mathfrak{q}$  is  a minimal primary decomposition of  
$I \cap \mathfrak{q}$.   Let  $I$ and $I \cap \mathfrak{q}$ be  normally torsion-free.  
Then, for all $m,n \geq 1$, we have $I^n(I\cap \mathfrak{q})^m=I^{n+m} \cap \mathfrak{q}^m$. 
  \end{lemma}

\begin{proof}
Fix $m,n \geq 1$. We first  verify that $I^n(I\cap \mathfrak{q})^m\subseteq I^{n+m} \cap \mathfrak{q}^m$.  According to \cite[Remarks 2.28(i)]{sharp}, we have 
the following inclusions 
$$I^n(I\cap \mathfrak{q})^m\subseteq I^{n} \cap (I \cap \mathfrak{q})^m \subseteq (I\cap \mathfrak{q})^m\subseteq  \mathfrak{q}^m.$$ 
Moreover,  since $(I \cap \mathfrak{q})^m \subseteq I^m$, we get $I^n(I\cap \mathfrak{q})^m\subseteq I^{n}.I^{m}=I^{m+n}$. This gives that  
$I^n(I\cap \mathfrak{q})^m\subseteq I^{n+m} \cap \mathfrak{q}^m$. 
  It thus  remains to prove the reverse inclusion. To accomplish  this, take an arbitrary monomial $u$ in $I^{n+m} \cap \mathfrak{q}^m$. 
  It follows from $u\in I^{n+m}$ that there exist  monomials $g_1, \ldots, g_n, h_1, \ldots, h_m\in \mathcal{G}(I)$ and some monomial $M$ such that 
  $u=g_1 \cdots g_n h_1 \cdots h_m M$. We also deduce from $u\in \mathfrak{q}^m$ that $u=x_{i_1} \cdots x_{i_m}w$, where 
 $x_{i_1},  \ldots, x_{i_m}\in \mathfrak{q}$ and  some monomial  $w$. Hence, we get   $x_{i_1} \cdots x_{i_m} \mid g_1 \cdots g_n h_1 \cdots h_m M$. 
 Without loss of generality, we partite the set $\{x_{i_1}, \ldots, x_{i_m}\}$ into $S_1, \ldots, S_a$, where $a\leq n$, 
 $S_{a+1}, \ldots, S_b$, where $b-a\leq m$, $S_{b+1}$ sets such that  $\bigcup_{i=1}^{b+1}S_i=\{x_{i_1}, \ldots, x_{i_m}\}$ and 
 for all $x_j \in S_k$ with $1\leq k \leq a$ (respectively, $x_j \in S_z$ with $a+1 \leq z \leq b$, and $x_j \in S_{b+1}$), we have 
 $x_j \mid g_k$ (respectively, $x_j \mid h_z$, and $x_j \mid M$). This implies that $g_k \in \mathfrak{q}^{|S_k|}$ (respectively, 
 $h_z \in \mathfrak{q}^{|S_z|}$, and $M\in \mathfrak{q}^{|S_{b+1}|}$). Now, we can distribute the rest  $n-a+m-b+a=m+n-b$ monomials of types $g_{i}$'s and  $h_j$'s such that, for $k=1, \ldots, a$,  $g_kH_k \in \mathfrak{q}^{|S_k|}\cap I^{|S_k|}$, where $H_k$ is the  product of $|S_k|-1$  monomials of  $m+n-b$ monomials of types $g_{i}$'s and  $h_j$'s (respectively, for $z=a+1, \ldots, b$, 
 $h_z T_z\in \mathfrak{q}^{|S_z|}\cap I^{|S_z|}$,  where $T_z$ is the  product of $|S_z|-1$  monomials of  $m+n-b$ monomials of types $g_{i}$'s and  $h_j$'s,    and $ML\in \mathfrak{q}^{|S_{b+1}|}\cap I^{|S_{b+1}|}$,  where $L$ is the  product of $|S_{b+1}|$  monomials of  $m+n-b$ monomials of types $g_{i}$'s and  $h_j$'s).
 Since  $I$ and  $I \cap \mathfrak{q}$ are   normally torsion-free, one can deduce from Facts \ref{proposition 4.3.29} and \ref{exercise 6.1.25} that, for all $\ell \geq 1$, we obtain $I^\ell=I^{(\ell)}=\bigcap_{\mathfrak{p}\in \mathrm{Min}(I)} \mathfrak{p}^\ell$ and  
  $(I\cap \mathfrak{q})^\ell=(I\cap \mathfrak{q})^{(\ell)}=\bigcap_{\mathfrak{p}\in \mathrm{Min}(I\cap \mathfrak{q})} \mathfrak{p}^\ell.$ 
By virtue of   $\bigcap_{\mathfrak{p}\in \mathrm{Ass}(I)}\mathfrak{p} \cap \mathfrak{q}$  is  a  minimal primary decomposition, we get  $ \mathrm{Min}(I\cap \mathfrak{q})=\mathrm{Min}(I) \cup \{\mathfrak{q}\}$, and so  $(I\cap \mathfrak{q})^\ell=I^\ell \cap \mathfrak{q}^\ell$.  
 In particular, we have  $g_kH_k \in (\mathfrak{q}\cap I)^{|S_k|}$ for $k=1, \ldots, a$  (respectively, 
 $h_zT_z \in (\mathfrak{q}\cap I)^{|S_z|}$   for $z=a+1, \ldots, b$, and $ML\in (\mathfrak{q}\cap I)^{|S_{b+1}|}$).
 It follows from $\bigcup_{i=1}^{b+1}S_i=\{x_{i_1}, \ldots, x_{i_m}\}$ that $\sum_{i=1}^{b+1}|S_i|=m$. This gives rise to 
 $u\in I^n(I\cap \mathfrak{q})^m$. 
Consequently,  for all $m,n \geq 1$, we obtain  $I^n(I\cap \mathfrak{q})^m=I^{n+m} \cap \mathfrak{q}^m$, as claimed.  
\end{proof}


\begin{proposition}\label{PRODUCT}
Let $I\subset R=K[x_1, \ldots, x_r]$  be a  square-free monomial ideal  and  $\mathfrak{q}$ be  a prime monomial ideal in $R$ such that $\bigcap_{\mathfrak{p}\in \mathrm{Ass}(I)}\mathfrak{p} \cap \mathfrak{q}$  is  a minimal primary decomposition of  
$I \cap \mathfrak{q}$.   Let   $I$ and $I \cap \mathfrak{q}$ be  normally torsion-free. Then, for all $m,n \geq 1$, the following statements hold:
 \begin{itemize}
 \item[(i)]     $\mathrm{Ass}(I^n(I\cap \mathfrak{q})^m)=\mathrm{Min}(I) \cup \{\mathfrak{q}\}$. 
 \item[(ii)] The monomial ideal  $I^n(I\cap \mathfrak{q})^m$ is normally torsion-free. 
  \end{itemize}
 \end{proposition}

\begin{proof}
(i) Fix $m,n \geq 1$.  It follows from  Lemma \ref{PRODUCT1} that  $I^n(I\cap \mathfrak{q})^m=I^{n+m} \cap \mathfrak{q}^m$.  Thus, we must show  that  $\mathrm{Ass}(I^{n+m} \cap \mathfrak{q}^m)=\mathrm{Min}(I) \cup \{\mathfrak{q}\}$. To do this, let $\mathrm{Ass}(I)=\{\mathfrak{p}_1, \mathfrak{p}_2, 
\ldots, \mathfrak{p}_t\}$. Since $I$ is normally torsion-free,  Facts \ref{proposition 4.3.29} and \ref{exercise 6.1.25} imply  
 $I^{n+m}= \mathfrak{p}^{n+m}_1 \cap \mathfrak{p}^{n+m}_2 \cap \cdots \cap \mathfrak{p}^{n+m}_t$ is a minimal primary decomposition of $I^{n+m}$. We show  $\mathfrak{p}^{n+m}_1 \cap \mathfrak{p}^{n+m}_2 \cap \cdots \cap \mathfrak{p}^{n+m}_t \cap \mathfrak{q}^m$ 
  is a minimal primary decomposition  of  $I^{n+m} \cap \mathfrak{q}^m$. For  this, we  prove 
    $\mathfrak{p}^{n+m}_1 \cap \mathfrak{p}^{n+m}_2 \cap \cdots \cap \mathfrak{p}^{n+m}_t \nsubseteq \mathfrak{q}^m$ and also  
  $\bigcap_{i=1, i\neq j}^{t} \mathfrak{p}^{n+m}_i \cap \mathfrak{q}^m \nsubseteq \mathfrak{p}^{n+m}_j,$ 
 where $j=1, \ldots, t$.     If $\mathfrak{p}^{n+m}_1 \cap \mathfrak{p}^{n+m}_2 \cap \cdots \cap \mathfrak{p}^{n+m}_t \subseteq \mathfrak{q}^m$, then $\sqrt{I^{n+m}} \subseteq \sqrt{\mathfrak{q}^m}$, and so $\sqrt{I}\subseteq \mathfrak{q}$. Since $I=\sqrt{I}$, we get $I\subseteq \mathfrak{q}$, which  contradicts the fact that $\bigcap_{\mathfrak{p}\in \mathrm{Ass}(I)}\mathfrak{p} \cap \mathfrak{q}$  is  a minimal primary decomposition. Suppose, on the contrary, that    $\bigcap_{i=1, i\neq j}^{t} \mathfrak{p}^{n+m}_i \cap \mathfrak{q}^m \subseteq \mathfrak{p}^{n+m}_j$ for  some $1\leq j \leq t$. As   $\mathfrak{p}^{n+m}_1 \cap \mathfrak{p}^{n+m}_2 \cap \cdots \cap \mathfrak{p}^{n+m}_t$ is a minimal primary decomposition of $I^{n+m}$, this implies that there exists a monomial  $u \in \bigcap_{i=1, i\neq j}^{t} \mathfrak{p}^{n+m}_i \setminus \mathfrak{p}^{n+m}_j$.
 On account of  $\bigcap_{\mathfrak{p}\in \mathrm{Ass}(I)}\mathfrak{p} \cap \mathfrak{q}$ is a minimal primary decomposition, we obtain   $\mathfrak{q}\nsubseteq \mathfrak{p}_j$. Hence, there exists some monomial  $v\in \mathfrak{q}\setminus  \mathfrak{p}_j$. This yields that  $uv^m \in \bigcap_{i=1, i\neq j}^{t} \mathfrak{p}^{n+m}_i  \cap \mathfrak{q}^m$, and so $uv^m \in \mathfrak{p}^{n+m}_j$. In view of  \cite[Corollary 6.1.8]{V1}, we know that $ \mathfrak{p}^{n+m}_j$ is  $\mathfrak{p}_j$-primary. We can now derive from $uv^m \in  \mathfrak{p}^{n+m}_j$ that $u\in  \mathfrak{p}^{n+m}_j$ or $v^m\in \sqrt{ \mathfrak{p}^{n+m}_j}$.  We therefore  get,  $u\in \mathfrak{p}^{n+m}_j$ or $v\in \mathfrak{p}_j$. This leads to a contradiction. Hence,  $\bigcap_{i=1}^{t}\mathfrak{p}^{n+m}_i  \cap \mathfrak{q}^m$   is a minimal primary decomposition  of  $I^{n+m} \cap \mathfrak{q}^m$.  Consequently,   $\mathrm{Ass}(I^n(I\cap \mathfrak{q})^m)=\mathrm{Min}(I) \cup \{\mathfrak{q}\}.$

(ii) It is well-known that if  $A$ and $B$ are two ideals in a commutative ring $S$, then by using induction on $\alpha, \beta \geq 1$, 
it can be shown that $A^{\alpha}B^{\beta}=B^{\beta}A^{\alpha}$ for all $\alpha, \beta \geq 1$. In particular, this fact implies that $(AB)^{\lambda}=A^{\lambda}B^{\lambda}$ for all $\lambda \geq 1$. Now, set $L:=I^n(I\cap \mathfrak{q})^m$ and fix $t\geq 1$. Since the polynomial ring $R=K[x_1, \ldots, x_n]$ is commutative (in particular,  due to \cite[Remarks 2.28(i)]{sharp}, we have $I(I\cap \mathfrak{q})=(I\cap \mathfrak{q})I$),
 this gives rise to the following equalities
$$L^t= (I^n(I\cap \mathfrak{q})^m)^t=(I^n)^t ((I\cap \mathfrak{q})^m)^t=I^{nt}(I\cap \mathfrak{q})^{mt}.$$
It follows now  from part (i) that  
$$\mathrm{Ass}(L^t)=\mathrm{Ass}(I^{nt}(I \cap \mathfrak{q})^{mt})=\mathrm{Min}(I) \cup \{\mathfrak{q}\}=\mathrm{Ass}(L).$$
This means that $L$ is normally torsion-free, as claimed. 
\end{proof}

To establish Lemma \ref{NTF1}, we need to use Lemma \ref{Ass-contraction}. For this purpose, we state and prove it in the next  lemma.

\begin{lemma} \label{Ass-contraction} 
Let $I\subset R=K[x_1, \ldots, x_n]$  be a  normally torsion-free square-free monomial ideal,  $\mathfrak{q}$ be  a prime monomial ideal in $R$,  and  $h$ be a square-free monomial in $R$ with  $\mathrm{supp}(h) \cap  (\mathrm{supp}(\mathfrak{q}) \cup \mathrm{supp}(I))=\emptyset$  such that $\bigcap_{\mathfrak{p}\in \mathrm{Ass}(I)}\mathfrak{p} \cap \bigcap_{x_r\in \mathrm{supp}(h)}(\mathfrak{q}, x_r)$  is  a minimal primary decomposition of  $L:=I \cap (\mathfrak{q}, h)$.   If  $\mathfrak{p}\in\mathrm{Ass}(L^s)\setminus \mathrm{Min}(L)$ for some $s\geq 2$, 
then  $(\mathfrak{q}, x_r) \subsetneq \mathfrak{p}$ for  some  $x_r \in \mathrm{supp}(h)$. 
\end{lemma}

\begin{proof}
Let $\mathfrak{p}\in\mathrm{Ass}(L^s)\setminus \mathrm{Min}(L)$ for some $s\geq 2$. Since $\mathfrak{p}\notin  \mathrm{Min}(L)$, this implies that $(\mathfrak{q}, x_r) \neq \mathfrak{p}$ for any $x_r \in \mathrm{supp}(h)$.  We claim that there exists some $x_r \in \mathrm{supp}(h)$ such that 
$x_r \in \mathfrak{p}$.  Suppose, on the contrary, that $x_r \notin \mathfrak{p}$ for any $x_r \in \mathrm{supp}(h)$. 
 Due to $\mathfrak{p}\in\mathrm{Ass}(L^s)$, it follows readily from   Fact \ref{lemma 4.6(viii)}  that $\mathfrak{p} \in \mathrm{Ass}(L(\mathfrak{p})^s)$.
Thanks to  $x_r \notin \mathfrak{p}$ for any $x_r \in \mathrm{supp}(h)$,  we get  $(\mathfrak{q}, h)(\mathfrak{p})=(1)$, and by virtue of 
 $L(\mathfrak{p})=I(\mathfrak{p}) \cap (\mathfrak{q}, h)(\mathfrak{p})$, this leads to $L(\mathfrak{p})=I(\mathfrak{p})$. In addition, the assumption $\mathrm{supp}(h) \cap  (\mathrm{supp}(\mathfrak{q}) \cup \mathrm{supp}(I))=\emptyset$ yields that  $L=I\cap \mathfrak{q} +hI$. It  follows  immediately  from $L \subseteq \mathfrak{p}$ and $x_r \notin \mathfrak{p}$ for any $x_r \in \mathrm{supp}(h)$ that    $I \subseteq \mathfrak{p}$. 
 One can now  deduce from  Fact \ref{theorem 3.15}   that $I(\mathfrak{p})$ is normally torsion-free, and so $L(\mathfrak{p})$ is normally torsion-free as well.  We  therefore have $\mathfrak{p} \in \mathrm{Ass}(L(\mathfrak{p}))$, and   Fact \ref{lemma 4.6(viii)} gives that  $\mathfrak{p} \in \mathrm{Ass}(L)$. This  implies  that $\mathfrak{p}\in \mathrm{Min}(L)$, which is a contradiction. Accordingly, there exists some $x_r \in \mathrm{supp}(h)$ such that 
$x_r \in \mathfrak{p}$, say $x_{r_0}$.  As  $L^s \subseteq \mathfrak{p}$  and $x_{r_0} \in \mathfrak{p}$, this forces us to obtain   $(\mathfrak{q}, x_{r_0}) \subsetneq \mathfrak{p}$,  as claimed.    
\end{proof}


We are in a position to provide  the first main result of this paper in the following lemma. Before stating it, we need to recall the following lemma. 

\begin{lemma}(\cite[Lemma 2.1]{KHN2})   \label{Kaplansky}
Let $S$ be a commutative ring and let $a_1,\ldots,a_m$ be elements constituting a permutable  $S$-sequence. Let $J$
 be an ideal generated by monomials in $a_{t+1},\ldots,a_m$  for some $t\in\mathbb{N}$ with $1\leq t \leq m-1$. Then $(J:_S a_1^{n_1} a_2^{n_2}
 \cdots a_t^{n_t})=J$ for all $n_1,n_2,\ldots,n_t\in\mathbb{N}$.
\end{lemma}

\begin{lemma} \label{NTF1}
Let $I\subset R=K[x_1, \ldots, x_n]$  be a  normally torsion-free square-free monomial ideal,  $\mathfrak{q}$ be  a prime monomial ideal in $R$,  and  $h$ be a square-free monomial in $R$ with  $\mathrm{supp}(h) \cap  (\mathrm{supp}(\mathfrak{q}) \cup \mathrm{supp}(I))=\emptyset$  
 such that  $\bigcap_{\mathfrak{p}\in \mathrm{Ass}(I)}\mathfrak{p} \cap \bigcap_{x_r\in \mathrm{supp}(h)}(\mathfrak{q}, x_r)$   (respectively,  $\bigcap_{\mathfrak{p}\in \mathrm{Ass}(I)}\mathfrak{p} \cap \mathfrak{q}$)  is a minimal primary decomposition of $I \cap (\mathfrak{q}, h)$ (respectively,  $I \cap \mathfrak{q}$).  Let   $I$ and $I \cap \mathfrak{q}$ be normally torsion-free. Then 
$L:=I \cap (\mathfrak{q}, h)$  is normally torsion-free.  
\end{lemma}

\begin{proof}
Let $\mathfrak{p}\in\mathrm{Ass}(L^s)\setminus \mathrm{Min}(L)$ for some $s\geq 2$. It follows at once from  Lemma  \ref{Ass-contraction} that 
$(\mathfrak{q}, x_r) \subsetneq \mathfrak{p}$   for  some  $x_r \in \mathrm{supp}(h)$.  Hence, we can view $\mathfrak{p}$ as 
 $\mathfrak{p}=(\mathfrak{q}', \mathfrak{q}'')$ with $\mathfrak{q} \subsetneq \mathfrak{q}'\subseteq \mathrm{supp}(\mathfrak{q}) \cup 
 \mathrm{supp}(I)$ and $\mathfrak{q}'' \subseteq \mathrm{supp}(h)$. 
Since  $L=I \cap (\mathfrak{q}, h)$ and $\mathrm{supp}(h) \cap  (\mathrm{supp}(\mathfrak{q}) \cup \mathrm{supp}(I))=\emptyset$, we get 
 $L=I \cap \mathfrak{q} + hI$.   We first  claim that, for all $1 \leq \ell \leq s$, we have 
$(L^s:h^{\ell})=I^{\ell-1}(I\cap \mathfrak{q})^{s-\ell+1}+I^{\ell}L^{s-\ell}.$  
 To achieve this, we proceed by induction on $\ell$.  Assume  $\ell=1$. Using the binomial theorem, we get 
 $$L^s=(I\cap \mathfrak{q})^s + \sum_{i=1}^{s} (hI)^{i}(I\cap \mathfrak{q})^{s-i} \text{ and }
  L^{s-1}=\sum_{i=0}^{s-1} (hI)^{i}(I\cap \mathfrak{q})^{s-1-i}.$$

  Moreover, since  $\mathrm{supp}(h) \cap  (\mathrm{supp}(\mathfrak{q}) \cup \mathrm{supp}(I))=\emptyset$, we deduce from Lemma \ref{Kaplansky} that $((I\cap \mathfrak{q})^s:h)=(I\cap \mathfrak{q})^s$. Also, in view of \cite[Corollary 2.2]{KHN2} (when $I=0$), we have 
  $((hI)^{i}(I\cap \mathfrak{q})^{s-i}: h)=I(hI)^{i-1}(I\cap \mathfrak{q})^{s-i}$ for all $1\leq i \leq s$.  This allows us to  conclude  that 
  $$\sum_{i=1}^{s} ((hI)^{i}(I\cap \mathfrak{q})^{s-i}: h)=I\sum_{i=1}^{s} (hI)^{i-1}(I\cap \mathfrak{q})^{s-i}=
  I\sum_{i=0}^{s-1} (hI)^{i}(I\cap \mathfrak{q})^{s-1-i}.$$
  
We therefore  obtain  the following equalities     
 \begin{align*}
(L^s: h)&=((I\cap \mathfrak{q})^s:h) + \sum_{i=1}^{s} ((hI)^{i}(I\cap \mathfrak{q})^{s-i}: h) \\
 &=(I\cap \mathfrak{q})^s + I\sum_{i=0}^{s-1} (hI)^{i}(I\cap \mathfrak{q})^{s-1-i} \\
 &=(I\cap \mathfrak{q})^s + IL^{s-1}.
 \end{align*}
Hence, the assertion is true for the case in which $\ell=1$. Now, suppose, inductively, that $\ell >1$  and that the result has been shown for all values less than  $\ell$. In particular,  we deduce from the inductive hypothesis that 
 \begin{equation}\label{eq33}
(L^s:h^{\ell-1})=I^{\ell-2}(I\cap \mathfrak{q})^{s-\ell+2}+I^{\ell-1}L^{s-\ell+1}.
\end{equation}
Note that  a similar argument yields the following equality 
 \begin{equation}\label{eq44}
 (L^{s-\ell+1}:h)=(I\cap \mathfrak{q})^{s-\ell+1} + IL^{s-\ell}.
 \end{equation}
By using   (\ref{eq33})  and (\ref{eq44}),  we get the following equalities 
\begin{align*}
(L^s:h^{\ell}) &=(I^{\ell-2}(I\cap \mathfrak{q})^{s-\ell+2}+I^{\ell-1}L^{s-\ell+1}:h)\\
 &=I^{\ell-2}(I\cap \mathfrak{q})^{s-\ell+2} + I^{\ell-1}((I\cap \mathfrak{q})^{s-\ell+1} + IL^{s-\ell})\\ 
 &=I^{\ell-2}(I\cap \mathfrak{q})^{s-\ell+2} +  I^{\ell-1}(I\cap \mathfrak{q})^{s-\ell+1} + I^{\ell}L^{s-\ell} \\ 
 &=I^{\ell-1}(I\cap \mathfrak{q})^{s-\ell+1}+I^{\ell}L^{s-\ell}.
 \end{align*}
 This completes the inductive step, and so  the claim has been shown  by induction. In particular,  $(L^s:h^s)=I^s$. Moreover, we have the following equalities 
$$((L^s:h^{\ell}),h)=I^{\ell-1}(I\cap \mathfrak{q})^{s-\ell+1}+I^{\ell}L^{s-\ell} + hR=I^{\ell}(I\cap \mathfrak{q})^{s-\ell} + hR,$$
 for all    $1 \leq \ell \leq s$.  In what follows, our aim is to verify this claim that, for all $1 \leq \ell \leq s$, we have 
  $\mathfrak{p}\in \mathrm{Ass}(L^s:h^{\ell}).$  To establish  the assertion for $\ell=1$, we  consider the following short exact sequence,  where 
    $\psi(\overline{t})= h\overline{t}$  and  $\phi(\overline{t}) =\overline{t}$,
 \begin{equation}\label{eq55}
0\longrightarrow R/(L^s:h)\stackrel{\psi}{\longrightarrow} R/L^s \stackrel{\phi}{\longrightarrow} R/(L^s,h) \longrightarrow  0.
\end{equation} 
 We thus  deduce  from   (\ref{eq55}) that   $\mathrm{Ass}(L^s) \subseteq  \mathrm{Ass}(L^s:h) \cup  \mathrm{Ass}(L^s,h).$  
  This leads to  $\mathfrak{p} \in \mathrm{Ass}(L^s:h)$ or $\mathfrak{p} \in  \mathrm{Ass}(L^s,h).$ 
 Let   $\mathfrak{p} \in  \mathrm{Ass}(L^s,h).$ Also,  by  the binomial theorem, we have 
  $L^s=(I\cap \mathfrak{q})^s + \sum_{i=1}^{s} (hI)^{i}(I\cap \mathfrak{q})^{s-i}$. By virtue of  $(hI)^{i}(I\cap \mathfrak{q})^{s-i} \subseteq hR$ for all 
  $1\leq i \leq s$, this implies that  $\sum_{i=1}^{s} (hI)^{i}(I\cap \mathfrak{q})^{s-i} \subseteq hR$. This gives rise to   $L^s +hR= (I \cap \mathfrak{q})^s + hR$. 
  It follows immediately   from Fact  \ref{theorem 3.3}  that 
 $\mathfrak{q}'\in \mathrm{Ass}(I \cap \mathfrak{q})^s.$  Due to $I \cap \mathfrak{q}$ is normally torsion-free, this gives that 
  $\mathfrak{q}'\in \mathrm{Min}(I \cap \mathfrak{q}).$  Because $\bigcap_{\mathfrak{p}\in \mathrm{Ass}(I)}\mathfrak{p} \cap \mathfrak{q}$ is a   minimal primary decomposition and $\mathfrak{q} \subsetneq \mathfrak{q}'$, we get a  contradiction. Hence,  we obtain   $\mathfrak{p} \notin  \mathrm{Ass}(L^s,h)$, and so 
  $\mathfrak{p} \in \mathrm{Ass}(L^s:h)$. Consequently,  the claim holds for $\ell=1$. 
  Now, suppose, inductively, that $\ell >1$  and that the claim  has been proven  for all values less than  $\ell$.   
   We now  consider the following short exact sequence,
\begin{equation}\label{eq66}
0\longrightarrow R/(L^s:h^{\ell})\stackrel{\psi}{\longrightarrow} R/(L^s:h^{\ell-1}) \stackrel{\phi}{\longrightarrow}
 R/((L^s:h^{\ell-1}),h) \longrightarrow  0.
\end{equation}
It follows  rapidly  from the inductive hypothesis that $\mathfrak{p}\in \mathrm{Ass}(L^s:h^{\ell-1}).$
By  (\ref{eq66}), one has   $\mathrm{Ass}(L^s:h^{\ell-1}) \subseteq  \mathrm{Ass}(L^s:h^{\ell}) \cup  
\mathrm{Ass}((L^s:h^{\ell-1}),h)$, and so   $\mathfrak{p}\in \mathrm{Ass}(L^s:h^{\ell})$ or $\mathfrak{p}\in \mathrm{Ass}((L^s:h^{\ell-1}),h)$.   Let $\mathfrak{p}\in \mathrm{Ass}((L^s:h^{\ell-1}),h)$. Thanks to $((L^s:h^{\ell-1}),h)=I^{\ell-1}(I \cap \mathfrak{q})^{s-\ell+1}+hR$, 
 Fact \ref{theorem 3.3}  implies that $\mathfrak{q}'\in \mathrm{Ass}(I^{\ell-1}(I \cap \mathfrak{q})^{s-\ell+1}).$ 
According to Proposition \ref{PRODUCT}, we can derive that  $\mathfrak{q}'\in \mathrm{Min}(I) \cup \{\mathfrak{q}\}$.  
Since $\bigcap_{\mathfrak{p}\in \mathrm{Ass}(I)}\mathfrak{p} \cap \mathfrak{q}$ is a   minimal primary decomposition    and $\mathfrak{q} \subsetneq \mathfrak{q}'$,  this leads to a contradiction. We thus get $\mathfrak{p}\in \mathrm{Ass}(L^s:h^{\ell})$. This completes the inductive step, and hence  the claim has been proven   by induction. In particular, we have $\mathfrak{p}\in \mathrm{Ass}(L^s:h^{s})$, and so $\mathfrak{p}\in \mathrm{Ass}(I^s)$. Thanks to $I$ is normally torsion-free, this implies that $\mathfrak{p}\in \mathrm{Min}(I)$, a contradiction. 
  Consequently, $L$ is normally torsion-free, as desired. 
\end{proof}


 The  theorem below  has been shown in \cite[Theorem 2.3]{NQ} by combinatorial tools, but,  as an application of Lemma \ref{NTF1}, we re-prove it.  
To formulate Theorem \ref{DI-TREES}, we recall here the required definitions. Let $G$ be a simple finite undirected graph. 
  For each vertex $v \in V(G)$, the {\em closed neighborhood} of $v$ in $G$ is defined as follows: 
\[
N_G[v]=\{u \in V(G) : \{u,v\} \in E(G)\} \cup \{v\}.
\]
 A subset $S\subseteq V(G)$ is called a  {\em dominating set} of $G$ if $S \cap N_G[v] \neq \emptyset$   for all $v\in V(G)$. A dominating set is called {\it minimal} if it does not properly contain any other dominating set of $G$. A minimum dominating set of $G$ is a minimal dominating set with the smallest size.  The {\em dominating number} of $G$, denoted by $\gamma(G)$, is  the size of its minimum dominating set, that is,
\[
\gamma(G)=\mathrm{min} \{|S|: S \text{ is a minimal dominating set
of } G\}.
\]
  The {\em closed neighborhood ideal} of $G$ has been introduced as 
 \[
NI(G)= (\prod_{x_j \in N[x_i]} x_j:  x_i \in V(G)),
\]
and   the {\em dominating ideal}  of $G$ is defined as
\[
DI(G)=(\prod_{x_i \in S}x_i: S \text{ is a minimal dominating set of } G).
\]

Set $\gamma'_{G} = \mathrm{max}\{|S| : S \text{ is a minimal dominating set of } G\}$. Then  we have the following fact:
   \begin{fact} \label{lemma  2.2}(\cite[Lemma  2.2]{SM})
Any minimal prime ideal of $NI(G)$ is of the form $\langle x_{j_1}, \ldots,  x_{j_r} \rangle$, where $\{x_{j_1}, \ldots,  x_{j_r}\}$ is a minimal dominating 
set of $G$. Therefore  $ht(NI(G)) = \gamma(G)$, $bight(NI(G)) = \gamma'_G$ and  $ht(NI(G)^\vee) = \min\{\deg_G(x_i) : x_i \in  V(G)\} + 1$.
\end{fact}
It follows from   Fact \ref{lemma  2.2}  that $DI(G)$ is the Alexander dual of $NI(G)$, that is, $DI(G)=NI(G)^{\vee}$.
 In order to find  more properties of  closed neighborhood ideals and dominating  ideals, the reader may refer to
   \cite{NBR, NQBM, SM}.

      Before stating the next theorem, we require to recall the following proposition.

   \begin{proposition} (\cite[Proposition 3.9]{KHN3}) \label{proposition 3.9}
Let $T$ be a tree which is not a path graph. Then $T$ has a 
vertex $u\in V(T)$ with $\mathrm{deg}_Tu\geq 3$, and $w_1, w_2 \in L_T$ with $w_1\neq w_2$, and $v_1, \ldots, v_r, v_{r+1}, \ldots, v_s \in V(T)$ (possibly $s=0$) with $\mathrm{deg}_T{v_i}=2$ for each $i=1, \ldots, s$ such that  $$Q: w_1, v_1, \ldots, v_r, u, v_{r+1}, \ldots, v_s, w_2$$
 is a maximal path in $T$.  
\end{proposition}


\begin{theorem}\label{DI-TREES}
The dominating  ideals of trees are normally torsion-free.
\end{theorem}

  \begin{proof}
  Suppose that  $T=(V(T), E(T))$ is  a tree.  We use  induction on $n:=|E(T)|$. By  Fact \ref{lemma 3.12}  and the  fact that every prime monomial ideal is normally torsion-free,  one can establish  the claim  for any  $n\leq 3$.   Now, suppose, inductively, that $n>3$ and that the result has been  shown  for all trees whose edge sets have cardinalities less than $n$. Assume that  $T=(V(T), E(T))$ is  a tree with $n=|E(T)|$. Using the fact that $T$ is a tree and $n>3$, due to 
   Proposition \ref{proposition 3.9},   there exists at least one leaf  such as $w$  such that $w$ has a neighbor such as $v$ with  either $\mathrm{deg}_T(v)=2$ or $v$ has another leaf  in its neighbor set.  
   Since  $w\in N_T(v)$,   we can view the tree $T$ as  the   union of a  tree such as $T_1$  and the edge $\{v,w\}$, where $v$ is the cut-point of this union. Note that, based on  the inductive hypothesis,  $DI(T_1)$ is normally torsion-free.   
Here, one may consider the following three cases:
\vspace{1mm}

\textbf{Case 1.}  $v$ has another leaf  in its neighbor set.  Hence,   $\prod_{j\in N_{T_1}[v]}x_j \notin \mathcal{G}(NI(T_1))$.  
On account  of  Fact \ref{lemma  2.2}  and  the fact that $v\in N_{T_1}[v]$, we obtain the following equalities
\begin{align*}
DI(T)=& DI(T_1) \cap (x_w,x_j~:~ j\in N_{T_1}[v])R \cap (x_w, x_v)R\\
=& DI(T_1) \cap (x_w, x_v)R,
\end{align*}
where $R=K[x_{\alpha}~:~ \alpha \in V(T)]$. Since $x_w \notin \mathrm{supp}(DI(T_1))$, we can derive from  Fact \ref{lemma 2.17} 
   that $DI(T)$ is normally torsion-free. 
\vspace{1mm}

\textbf{Case 2.}   $\mathrm{deg}_T(v)=2$ and   $\prod_{j\in N_{T}[c]}x_j \notin \mathcal{G}(NI(T))$, where  $N_{T}(v)=\{c\}$. 
This means that there is a leaf in the neighbor set of $c$.    Let $L$ be a square-free monomial ideal such that  $\mathcal{G}(L)= \mathcal{G}(NI(T_1))\setminus \{x_cx_v \}$.   We thus get   
$$DI(T)=L^\vee \cap (x_w,x_c,x_v)R \cap (x_w,x_v)R=L^\vee \cap (x_w,x_v)R,$$  where $R=K[x_{\alpha}~:~ \alpha \in V(T)]$ and 
  $L^{\vee}$   denotes  the Alexander dual of $L$. Thanks to $L^{\vee}=DI(T_1)/x_v$ and $DI(T_1)$ is normally torsion-free, it follows at once from 
  Fact \ref{theorem 3.19}   that $L^{\vee}$ is normally torsion-free. As   $x_v, x_w\notin \mathrm{supp}(L^{\vee})$, Fact \ref{proposition 3.3}
       gives that    $DI(T)$ is normally torsion-free.
 \vspace{1mm}

\textbf{Case 3.}   $\mathrm{deg}_T(v)=2$ and  $\prod_{j\in N_{T}[c]}x_j \in \mathcal{G}(NI(T))$, where  $N_{T_1}(v)=\{c\}$.  
Notice that it follows readily  from $N_{T_1}(v)=\{c\}$ that one  can view the tree $T_1$ as  the   union of a  tree such as $T_2$  and the edge $\{c,v\}$, where $c$ is the cut-point of this union.
 Assume   $L_1$ and $L_2$ are  two  square-free monomial ideals such that 
 $\mathcal{G}(L_1)=\mathcal{G}(NI(T_1))\setminus \{x_cx_v \}\; \; \text{and}\; \; \mathcal{G}(L_2)=\mathcal{G}(L_1)\cup \{\prod_{j\in N_{T_1}[c]}x_j\}.$
  In view of  Fact \ref{lemma  2.2},  one gains  the following equalities  
 $$DI(T_1)=L^{\vee}_1\cap (x_c, x_v)\; \; \text{and}\;\; L^{\vee}_2=L^{\vee}_1 \cap (x_j~:~ j\in N_{T_1}[c]),$$
 where $L^{\vee}_1$ (respectively, $L^{\vee}_2$)  denotes  the Alexander dual of $L_1$ (respectively, $L_2$).
 It is not hard to check that $L^{\vee}_1=DI(T_1)/x_v$. Due to   $DI(T_1)$ is normally torsion-free, one can deduce from Fact \ref{theorem 3.19}  that 
 $L^{\vee}_1$ is normally torsion-free. Put $\mathfrak{q}:=(x_j~:~ j\in N_{T_2}[c])$. Then we have 
 $L^{\vee}_2=L^{\vee}_1 \cap (\mathfrak{q}, x_v).$
 On account of  $DI(T_2)$ is normally torsion-free and $DI(T_2)=L^{\vee}_1\cap \mathfrak{q}$, this implies that $L^{\vee}_1\cap \mathfrak{q}$
 is normally torsion-free.  Notice that $x_v \notin \mathrm{supp}(\mathfrak{q}) \cup \mathrm{supp}(L^{\vee}_1)$. 
 By virtue of Lemma \ref{NTF1}, one can immediately deduce that    $L^{\vee}_2$ is normally torsion-free.   Once again, Fact \ref{lemma  2.2}  yields that 
$$DI(T)=L^{\vee}_2\cap (x_w, x_j~:~ j\in N_{T_1}[v])R \cap (x_w,x_v)R=L^{\vee}_2 \cap (x_w,x_v)R,$$  where $R=K[x_{\alpha}~:~ \alpha \in V(T)]$.  
Due to  $x_w \notin \mathrm{supp}(L^{\vee}_2)$ and $L^{\vee}_2$ is normally torsion-free,  one  can immediately conclude from 
 Fact \ref{lemma 2.17}    that   $DI(T)$ is normally torsion-free.
This completes the inductive step, and hence  the claim has been proven by induction. 
\end{proof}

\section{Some results on the embedded  associated primes} 

This section  has two folds. We first give some results on the embedded  associated prime ideals of powers of square-free monomial ideals. Next, we turn our attention to study when a square-free monomial ideal is minimally not normally torsion-free. Particularly, we introduce  a class of square-free monomial ideals, which are minimally not normally torsion-free.

To understand   Lemma    \ref{Ass-maximal1} and     Proposition \ref{MNNT},  one  needs  to recollect  some definitions. Assume that $I$ is  a square-free monomial ideal and $\Gamma \subseteq \mathcal{G}(I)$, where $\mathcal{G}(I)$ denotes the unique minimal set of monomial generators of the  monomial ideal $I$. We say that $\Gamma$ is an {\it independent} set in $I$ if $\mathrm{gcd}(f,g)=1$ for each $f,g\in \Gamma$ with $f\neq g$. We denote the maximum cardinality of an independent set in $I$ by $\beta_1(I)$. 
Furthermore, for $1\leq j \leq n$,   the { \it contraction} of $I$ at $x_j$, denoted by $I/x_j$,  is obtained by setting $x_j=1$ in  each $u\in \mathcal{G}(I)$.  In addition, the  { \it deletion} of $I$ at $x_j$, denoted by $I\setminus x_j$,     is formed  by putting  $x_j=0$ in  each $u\in \mathcal{G}(I)$, refer to   \cite[page 303]{HM}.\par

We are ready to give the second main result of this paper in the following lemma, which  is an affirmative answer to \cite[Question 3.11]{NQ1}. Before stating this lemma, we require to recall the following proposition. 

\begin{proposition}\label{proposition 4.8} (\cite[Proposition 4.8]{SNQ})
Let $I$ be a square-free monomial ideal in a polynomial ring $R=K[x_1, \ldots, x_n]$ over a field $K$.  Let   
$\mathfrak{p}=(I^s:h)\in \mathrm{Ass}(R/I^s)$ for some positive integer $s$ and some monomial $h$ in $R$. Then $\mathrm{deg}_{x_i}h\leq s-1$ for each $i=1, \ldots, n$. 
\end{proposition}

\begin{lemma} \label{Ass-maximal1}
Let $I \subset R=K[x_1, \ldots, x_n]$ be a  square-free monomial ideal  and   $\mathfrak{q}$ be  a prime monomial ideal in $R$. Let  $L:=I \cap (\mathfrak{q}, x_r)$  such that   $x_r\notin \mathrm{supp}(\mathfrak{q}) \cup \mathrm{supp}(I)$ and  $\mathrm{Ass}(L)=\mathrm{Ass}(I) \cup \{(\mathfrak{q}, x_r)\}$.   If    $\mathfrak{p}=(L^s:x_r^{\theta}v)$ for some $s\geq 2$,   $\theta \geq 0$, and monomial $v$  with 
$x_r\nmid v$, then  $\mathfrak{p}\setminus  x_r \in  \mathrm{Ass}(I^{\theta}(I\cap \mathfrak{q})^{s-\theta})$. 
\end{lemma}

\begin{proof}
Let  $\mathfrak{p}=(L^s:x_r^{\theta}v)$ for some $s\geq 2$,   $\theta \geq 0$, and monomial $v$  with 
$x_r\nmid v$. 
   Since $x_r^{\theta}v\notin L^s$ and $L^s=\sum_{i=0}^{s} (x_rI)^{i}(I\cap \mathfrak{q})^{s-i}$, this gives  that $x_r^{\theta}v \notin (x_rI)^{i} (I\cap \mathfrak{q})^{s-i}$ for all $i=0, \ldots, s$. Thanks to $L$ is a square-free monomial ideal, one can deduce from Proposition \ref{proposition 4.8}
    that $\theta  \leq s-1$. In particular, we have $x_r^{\theta}v \notin (x_rI)^{\theta}(I\cap \mathfrak{q})^{s-\theta}$, and so 
  $v \notin I^{\theta}(I\cap \mathfrak{q})^{s-\theta}$. Now, let $x_t$ be an arbitrary element in $\mathfrak{p}\setminus  x_r$. It follows from 
  $\mathfrak{p}=(L^s: x_r^{\theta}v)$ that $x_tx_r^{\theta}v\in L^s$, and so $x_tv\in (L^s:x_r^{\theta})$. Because $x_r \nmid x_tv$, this yields  that 
  $x_tv\in (L^s:x_r^{\theta})\setminus x_r$. By remembering the following chain
  \begin{equation}\label{eq77}
(I\cap \mathfrak{q})^s \subseteq I(I\cap \mathfrak{q})^{s-1} \subseteq I^2(I\cap \mathfrak{q})^{s-2} \subseteq \cdots \subseteq I^{s-1}(I\cap \mathfrak{q}) \subseteq I^s, 
\end{equation}
we get $x_tv\in I^{\theta}(I\cap \mathfrak{q})^{s-\theta}$. We can conclude from $v\notin I^{\theta}(I\cap \mathfrak{q})^{s-\theta}$ and 
$x_tv \in I^{\theta}(I\cap \mathfrak{q})^{s-\theta}$ for all $x_t \in \mathfrak{p}\setminus  x_r$  that 
$(I^{\theta}(I\cap \mathfrak{q})^{s-\theta}:v)=\mathfrak{p}\setminus  x_r$.  This gives rise to  $\mathfrak{p}\setminus  x_r  \in  \mathrm{Ass}(I^{\theta}(I\cap \mathfrak{q})^{s-\theta})$, and the proof is complete. 
 \end{proof}


As an immediate consequence of Lemma \ref{Ass-maximal1}, we get  the  corollary below. 
\begin{corollary} \label{Cor.Ass-maximal1}
Let $I \subset R=K[x_1, \ldots, x_n]$ be a  square-free monomial ideal  and   $\mathfrak{q}$ be  a prime monomial ideal in $R$. Let  $L:=I \cap (\mathfrak{q}, x_r)$  such that   $x_r\notin \mathrm{supp}(\mathfrak{q}) \cup \mathrm{supp}(I)$ and  $\mathrm{Ass}(L)=\mathrm{Ass}(I) \cup \{(\mathfrak{q}, x_r)\}$.   If    $\mathfrak{p}=(x_1, \ldots, x_n)=(L^s:v)$ for some $s\geq 2$   with  
$x_r\nmid v$, then   $\mathfrak{p}\setminus  x_r \in  \mathrm{Ass}((I\cap \mathfrak{q})^s)$. 
\end{corollary}


\begin{example}\label{Exam.1}
Let   $R=K[x_1, \ldots, x_8]$, $\mathfrak{q}:=(x_1,x_3,x_7)$, $x_r:=x_8$, and  
\begin{align*}
I =&(x_1,x_4,x_7)\cap (x_1,x_4,x_2)\cap (x_2,x_4,x_6)\cap (x_2,x_6,x_5)\cap (x_3,x_5,x_6)\\
&\cap (x_5,x_7,x_3)\\
=&(x_1x_2x_3,x_2x_3x_4,x_1x_2x_5,x_4x_5,x_1x_3x_6,x_3x_4x_6,x_1x_5x_6,x_2x_3x_7,\\
& x_2x_5x_7,x_1x_6x_7,  x_2x_6x_7,x_4x_6x_7).
\end{align*}
Let $L:=I\cap (\mathfrak{q}, x_8)$. Using {\it Macaulay2} \cite{GS}, we obtain $(L^2: v)=(x_1, \ldots, x_8)$, where 
$v=x_2x_3x_4x_5x_7$. Since $x_8\nmid v$, it  follows immediately  from Corollary \ref{Cor.Ass-maximal1} that 
$(x_1, \ldots, x_7)\in \mathrm{Ass}((I\cap \mathfrak{q})^2)$.
\end{example}

In the  following result   we are going to introduce a class of square-free monomial ideals which are minimally not normally torsion-free. 
An ideal $I'$ is  a {\it minor}   of  a square-free monomial ideal $I$ if $I'$ can be obtained from $I$ by a sequence of deletions and contractions in any order, see \cite[Definition 6.5.3]{V1}. Recall that   a square-free monomial ideal $I \subset R=K[x_1, \ldots, x_n]$ is said to be  {\it  minimally not normally torsion-free} when $I$  is not normally torsion-free but all its proper minors are  normally torsion-free, refer to \cite{HM}.  
Before stating the next proposition, we need to recall the following fact. 

\begin{fact}\label{lemma 4.12(vi)}(\cite[Lemma 4.12(vi)]{RNA})
Let  $I$ be a monomial ideal in a polynomial ring   $R=K[x_1,\ldots, x_n]$ over a field $K$,  and $1\leq j \leq n$.
 If $I=Q_1\cap \cdots \cap Q_s$ is a minimal primary decomposition of $I$ in $R$, then  $I\setminus x_j=(Q_1\setminus x_j)\cap \cdots \cap (Q_s\setminus x_j)$
  is a  primary decomposition of $I\setminus x_j$ in $R\setminus x_j$.
\end{fact}


\begin{proposition}\label{MNNT}
 Let $I\subset R=K[x_1, \ldots, x_n]$  be a square-free monomial ideal and  $\mathfrak{m}=(x_1, \ldots, x_n)$  
 such that   there exists  some $ m\geq 1$ with   $\mathrm{Ass}(I^s)=\mathrm{Min}(I)$ for all $1\leq s \leq m$, and  $\mathrm{Ass}(I^s)= \mathrm{Min}(I) {\cup} \{\mathfrak{m}\}$   for all $s\geq m+1$  such that  $\mathfrak{m}\setminus x_j \notin  I^s\setminus x_j$ for all $j=1, \ldots, n$. 
    Then $I$ is minimally not normally torsion-free. In particular, we have $m\geq \beta_1(I)$, where $\beta_1(I)$ denotes the maximum cardinality of an independent set in $I$. 
\end{proposition}

\begin{proof}
Since $\mathfrak{m}\in \mathrm{Ass}(I^s)$ for all $s\geq m+1$, this implies that $I$ is not normally torsion-free. To complete the proof, we have to show that $I/x_j$   and  $I\setminus x_j$  are normally torsion-free for all $1\leq j \leq n$. To do this, fix $1\leq j \leq n$.  We first verify that $I/x_j$ is normally torsion-free. To see this, we  need to show that  $\mathrm{Ass}((I/x_j)^k)  \subseteq \mathrm{Ass}(I/x_j)$  for all $k\geq 1$. Fix $k \geq 1$. 
By  Fact  \ref{lemma 4.10(viii)}, we have 
$$\mathrm{Ass}((I/x_j)^k)=\mathrm{Ass}(I^k/x_j)=\{\mathfrak{q} : \mathfrak{q} \in \mathrm{Ass}_{R}(R/I^k)~\mathrm{and}~ x_j\notin\mathfrak{q}\}.$$ Choose an arbitrary element $\mathfrak{q}\in \mathrm{Ass}((I/x_j)^k)$. We thus get  
$\mathfrak{q} \in \mathrm{Ass}_{R}(R/I^k)$ and $x_j\notin\mathfrak{q}$. Due to $x_j\notin\mathfrak{q}$, this yields that $\mathfrak{q}\neq \mathfrak{m}$, and so $\mathfrak{q}\in \mathrm{Min}(I)$. Accordingly, we have  $\mathfrak{q}\in   \mathrm{Ass}(I)$, and therefore 
   $\mathfrak{q}\in  \mathrm{Ass}(I/x_j)$. This gives rise to $I/x_j$ is  normally torsion-free. 
To finish the proof, we demonstrate that $I\setminus x_j$ is normally torsion-free. To accomplish this, we must establish 
$\mathrm{Ass}((I\setminus x_j)^k)  \subseteq \mathrm{Ass}(I \setminus x_j)$  for all $k\geq 1$. Fix $k \geq 1$. 
One can derive from  \cite[Definition 4.3.22]{V1} and  Fact \ref{exercise 6.1.25}  that $I^{(k)}=\bigcap_{\mathfrak{p}\in \mathrm{Ass}(I)}\mathfrak{p}^k$, and so 
$I^k=\bigcap_{\mathfrak{p}\in \mathrm{Ass}(I)}\mathfrak{p}^k$ when $1\leq k\leq m$, and 
$I^k=\bigcap_{\mathfrak{p}\in \mathrm{Ass}(I)}\mathfrak{p}^k \cap Q$ with $\sqrt{Q}=\mathfrak{m}$ when $k\geq m+1$. Let 
$\mathrm{Ass}(I)=\{\mathfrak{p}_1, \ldots, \mathfrak{p}_z\}$.  
From   Fact   \ref{lemma 4.12(vi)}, we have   a primary decomposition of $I^k\setminus x_j$  as follows:
$$I^k\setminus x_j  =  \bigcap_{i=1}^z(\mathfrak{p}_i\setminus x_j)^k \; \; \text{or} \; \;  I^k\setminus x_j = \bigcap_{i=1}^z(\mathfrak{p}_i\setminus x_j)^k \cap (Q\setminus x_j).$$
Let   $\mathrm{Ass}(I\setminus x_j)=\mathrm{Min}(I\setminus x_j)=\{{\mathfrak{p}}_{i_1}\setminus x_j, \ldots, {\mathfrak{p}}_{i_h}\setminus x_j\},$  where $\{i_1, \ldots, i_h\}\subseteq \{1, \ldots, z\}$. 
We note  if   $I^k\setminus x_j = \bigcap_{i=1}^z(\mathfrak{p}_i\setminus x_j)^k \cap (Q\setminus x_j),$ 
then   we have  $\bigcap_{t=1}^h(\mathfrak{p}_{i_t}\setminus x_j)^k  \subseteq  Q\setminus x_j$. 
Indeed, if it  were  $\bigcap_{t=1}^h(\mathfrak{p}_{i_t}\setminus x_j)^k  \nsubseteq  Q\setminus x_j$, then we can derive 
 that  $\mathfrak{m}\setminus x_j \in \mathrm{Ass}(I^k\setminus x_j)$. 
 This contradicts our assumption. We thus  get     $(I\setminus x_j)^k=({\mathfrak{p}}_{i_1}\setminus x_j)^k \cap \cdots \cap ({\mathfrak{p}}_{i_h}\setminus x_j)^k$   is  a   minimal primary decomposition of $(I\setminus x_j)^k$.
 As    $\mathrm{Min}(I\setminus x_j)=\{{\mathfrak{p}}_{i_1}\setminus x_j, \ldots, {\mathfrak{p}}_{i_h}\setminus x_j\},$ it follows from 
Fact \ref{exercise 6.1.25}  that $(I\setminus x_j)^{(k)}=({\mathfrak{p}}_{i_1}\setminus x_j)^k \cap \cdots \cap ({\mathfrak{p}}_{i_h}\setminus x_j)^k.$ 
This leads to  $(I\setminus x_j)^{(k)}=(I\setminus x_j)^k$. We can finally deduce from  Fact \ref{proposition 4.3.29}  that 
$I\setminus x_j$ is normally torsion-free. The last assertion can be deduced  from  Fact \ref{corollary 3.5}, and the proof is over. 
\end{proof}


We close this paper with the following example. In fact, by employing   Proposition   \ref{MNNT},  we are going to introduce a class of square-free monomial ideals, which are minimally not normally torsion-free. Before stating this example, we need to recall the following proposition. 

\begin{proposition}\label{proposition 3.6} (\cite[Proposition 3.6]{NKA})
 Suppose that $C_{2n+1}$ is  a cycle graph on the  vertex set $[2n+1]$, $R=K[x_1, \ldots, x_{2n+1}]$ is a  polynomial ring over a field $K$,
 and $\mathfrak{m}$ is the unique homogeneous maximal ideal  of $R$. Then  $\mathrm{Ass}_R(R/(J(C_{2n+1}))^s)= \mathrm{Ass}_R(R/J(C_{2n+1}))\cup \{\mathfrak{m}\},$  for all $s\geq 2$. In particular, 
 $$\mathrm{Ass}^\infty(J(C_{2n+1}))=\{(x_i, x_{i+1})~: ~ i=1, \ldots 2n\}\cup\{(x_{2n+1}, x_1)\}\cup \{\mathfrak{m}\}.$$
 \end{proposition}

\begin{example} \label{Exam.MNNT}
Let $C_{2n+1}=(V(C_{2n+1}), E(C_{2n+1}))$ be the odd cycle graph with  $V(C_{2n+1})=\{x_1, \ldots, x_{2n+1}\}$ and  
$$E(C_{2n+1})=\{\{x_i, x_{i+1}\} : i=1, \ldots, 2n\} \cup \{\{x_1, x_{2n+1}\}\}.$$  Also, let $J(C_{2n+1})$ denote   the cover ideal of $C_{2n+1}$ in 
the polynomial ring  $R=K[x_1, \ldots, x_{2n+1}]$ over a field $K$  and $\mathfrak{m}=(x_1, \ldots, x_{2n+1})$ be the maximal ideal. By 
Proposition \ref{proposition 3.6}, we know that 
$\mathrm{Ass}(J(C_{2n+1}))=\mathrm{Min}(J(C_{2n+1}))$ and $\mathrm{Ass}(J(C_{2n+1})^s)=\mathrm{Min}(J(C_{2n+1}))\cup \{\mathfrak{m}\}$ for all $s\geq 2$. It is easy  to see that for each $1\leq j \leq 2n+1$ with $x_0=x_{2n+1}$ and $x_{2n+2}=x_1$, we have 
  $J(C_{2n+1})\setminus x_j=(x_{j-1}) \cap (x_{j+1}) \cap J(P)$, where $P$ is  a path with $x_{j-1}, x_{j+1}\notin V(P)$. From   
 \cite[Corollary 2.6]{GRV}, $J(P)$ is normally torsion-free, and in view of  Fact \ref{proposition 3.3}, we dedude that $J(C_{2n+1})\setminus x_j$ is normally torsion-free, and so $\mathfrak{m}\setminus x_j \notin \mathrm{Ass}(J(C_{2n+1})^s\setminus x_j)$ for all $s$. It follows now from Proposition  \ref{MNNT} that $J(C_{2n+1})$ is minimally not normally torsion-free. In particular, we have $\beta_1(J(C_{2n+1}))=1$.  
\end{example}


\section*{ORCID}

Mehrdad ~Nasernejad: \textsc{https://orcid.org/0000-0003-1073-1934}\\
Veronica Crispin Qui$\mathrm{\tilde{n}}$onez:  \textsc{https://orcid.org/0000-0002-2146-7051}\\  
  Winfried  Hochstättler:  \textsc{https://orcid.org/0000-0001-7344-7143} \\

 
\section*{The conflict of interest and data availability statement}

We hereby declare  that this  manuscript has no associated  data and also there is no conflict of interest in this manuscript.


\section*{Acknowledgment.} 
First, the authors  are  deeply grateful to the anonymous referee for careful reading of the manuscript, and for  his/her  valuable suggestions which led to 
several  improvements in the quality of this paper. 
In addition, this paper was prepared when the first author  visited the Department of Mathematics of 
Uppsala University in 2023; in particular, he  would like to thank Uppsala 
University for its hospitality. Also, the second author was partially supported by the Lundstr\"om-\AA{m}an Foundation.



\end{document}